\documentclass[11pt]{article}

\usepackage{amsmath, amsfonts, amssymb, amsthm, graphicx,  color, enumitem}
\usepackage{fullpage}
\usepackage{eepic} 

\newcommand{\Prob}{\mathbb{P}}

\newcommand{\R}{\mathbb{R}}

\DeclareMathOperator{\Ai}{Ai}

\DeclareMathOperator{\Geom}{Geom}

\newtheorem{thm}{Theorem}[section]

\theoremstyle{definition}
\newtheorem{defn}{Definition}[section]

\theoremstyle{remark}
\newtheorem{rmk}{Remark}[section]

\newcommand{\Aip}{\mathcal{A}}
\newcommand{\Aipp}{\hat{\Aip}}

\newcommand\psymmU{%
\begin{picture}(1,1)(0,0)%
\allinethickness{0.5pt}%
\path(0,0)(0,1)(1,1)(1,0)(0,0)%
\end{picture}}
\newcommand\psymmUU{%
\begin{picture}(1,1)(0,0)%
\allinethickness{0.5pt}%
\path(0,0)(0,1)(1,1)(1,0)(0,0)%
\put(0.5,0.5){\makebox(0,0){$\cdot$}}%
\end{picture}}
\newcommand\psymmO{%
\begin{picture}(1,1)(0,0)%
\allinethickness{0.5pt}%
\path(0,0)(0,1)(1,1)(1,0)(0,0)%
\path(0,0)(1,1)%
\end{picture}}
\newcommand\psymmS{%
\begin{picture}(1,1)(0,0)%
\allinethickness{0.5pt}%
\path(0,0)(0,1)(1,1)(1,0)(0,0)%
\path(1,0)(0,1)%
\end{picture}}
\newcommand\psymmu{%
\begin{picture}(1,1)(0,0)%
\allinethickness{0.5pt}%
\path(0,0)(0,1)(1,1)(1,0)(0,0)%
\path(0,0)(1,1)%
\path(0,1)(1,0)%
\end{picture}}

\newbox\tsymmUbox
\newbox\tsymmUUbox
\newbox\tsymmObox
\newbox\tsymmSbox
\newbox\tsymmubox
\setbox\tsymmUbox =\hbox{\kern0.75pt\setlength{\unitlength}{6pt}\psymmU \kern0.75pt}
\setbox\tsymmUUbox=\hbox{\kern0.75pt\setlength{\unitlength}{6pt}\psymmUU\kern0.75pt}
\setbox\tsymmObox =\hbox{\kern0.75pt\setlength{\unitlength}{6pt}\psymmO \kern0.75pt}
\setbox\tsymmSbox =\hbox{\kern0.75pt\setlength{\unitlength}{6pt}\psymmS \kern0.75pt}
\setbox\tsymmubox =\hbox{\kern0.75pt\setlength{\unitlength}{6pt}\psymmu \kern0.75pt}

\newbox\symmUbox
\newbox\symmUUbox
\newbox\symmObox
\newbox\symmSbox
\newbox\symmubox
\setbox\symmUbox =\hbox{\kern0.75pt\setlength{\unitlength}{4.5pt}\psymmU \kern0.75pt}
\setbox\symmUUbox=\hbox{\kern0.75pt\setlength{\unitlength}{4.5pt}\psymmUU\kern0.75pt}
\setbox\symmObox =\hbox{\kern0.75pt\setlength{\unitlength}{4.5pt}\psymmO \kern0.75pt}
\setbox\symmSbox =\hbox{\kern0.75pt\setlength{\unitlength}{4.5pt}\psymmS \kern0.75pt}
\setbox\symmubox =\hbox{\kern0.75pt\setlength{\unitlength}{4.5pt}\psymmu \kern0.75pt}
\def\symmU{{\copy\symmUbox}}
\def\symmUU{{\copy\symmUUbox}}
\def\symmO{{\copy\symmObox}}
\def\symmS{{\copy\symmSbox}}
\def\symmu{{\copy\symmubox}}

\title{On the average of the Airy process and its time reversal}

\author{Jinho Baik\footnote{Department of Mathematics, University of Michigan,
Ann Arbor, MI, 48109, USA \newline email: \texttt{baik@umich.edu}}
and Zhipeng Liu\footnote{Department of Mathematics, University of
Michigan, Ann Arbor, MI, 48109, USA
\newline email: \texttt{zhpliu@umich.edu}}}

\date{\today}

\begin{document}

\maketitle

\begin{abstract} 
We show that the supremum of the average of the Airy process and its time reversal minus a parabola 
is distributed as the maximum of two independent GUE Tracy-Widom random variables. 
The proof is obtained by considering a directed last passage percolation model with a rotational symmetry in two different ways. 
We also review other known identities between the Airy process and the Tracy-Widom distributions.
\end{abstract}

\maketitle

\section{Introduction and result} \label{section:introduction}

\subsection{Airy process} 

We first introduce some notations. 

\begin{defn}\label{defn:TW}
Let $F_\beta(x)$, $\beta=1,2$ and $4$, denote the GOE, GUE and GSE Tracy-Widom distribution functions 
defined in \cite{Tracy-Widom94, Tracy-Widom96},  respectively.
\end{defn}

An explicit formula is included in Section~\ref{sec:formula} below. 

\begin{defn}
Let $\Aip(\tau)$, $\tau\in \R$, denote the Airy process.\footnote{It is also often called Airy$_2$ process and is denoted by $\Aip_2(\tau)$ in order to emphasize that its marginal is distributed as $F_2$. 
}
Set 
\begin{equation}
	\Aipp(\tau) := \Aip(\tau)-\tau^2.
\end{equation}
\end{defn}

The Airy process  was introduced in \cite{Prahofer-Spohn02}. 
It is believed that $\Aip(\tau)$ is an universal limit of the spatial fluctuations for models in the so-called KPZ (Kardar-Parisi-Zhang) class. 
The limit theorem to the Airy process is established for several special cases in $2$-dimensional directed last passage percolation, $1+1$ dimensional random growth, non-intersecting processes, and random matrices. 
See,  for example, \cite{Corwin12} and the references therein. 

The basic connection between the Airy process and the Tracy-Widom distribution is that the marginal distribution of $\Aip(\tau)$ at a fixed $\tau$ is $F_2(x)$. 
The joint distribution at finitely many times is also explicit and is given by a determinantal formula involving the Airy function. 
The Airy process is stationary but is not Markovian. 

In addition to the above basic connection, there are interesting identities between the supremum of a function of the Airy process and the Tracy-Widom distribution functions. 
The purpose of this paper is to establish one more such an identity. 
We first review two known identities. 
 
\begin{thm}[\cite{Johansson03}]\label{thm:A}
For every $x\in \R$, 
\begin{equation}\label{eq:thmA}
	\Prob\left[2^{2/3} \sup_{\tau\in \R}\Aipp(\tau)  \le x \right] = F_1(x).
\end{equation}
\end{thm}

\begin{thm}[\cite{Baik-Liu02}]\label{thm:Aab}
Let $\Aip^{(1)}$ and $\Aip^{(2)}$ be two independent Airy processes. 
Then 
\begin{equation}
	\Prob\left[  \sup_{\tau\in \R} 
	\left( \frac{\alpha^{1/3}\Aipp^{(1)}(\alpha^{-2/3}\tau)+\beta^{1/3}\Aipp^{(2)}(\beta^{-2/3}\tau)}{(\alpha+\beta)^{1/3}} \right)\le x \right] = F_2(x)
\end{equation}
for every $\alpha, \beta>0$ and for every $x\in \R$.
\end{thm}

The main result of this paper is the following identity. 

\begin{thm}\label{thm:AA}
For every $x\in \R$, 
\begin{equation}\label{eq:AA}
	\Prob\left[  \sup_{\tau\in \R} \left(\frac{\Aipp(\tau)+\Aipp(-\tau)}2\right) \le x \right] = F_2(x)^2.
\end{equation}
\end{thm}

Hence the supremum of the average of the Airy process and its time reversal minus a parabola 
is distributed as the maximum of two independent GUE Tracy-Widom random variables. 

Theorem~\ref{thm:A} was proved by Johansson in 2003. 
It was obtained indirectly by interpreting the so-called point-to-line last passage time of a directed last passage percolation (DLPP) model in two different ways.
Theorem~\ref{thm:Aab} was proved similarly by considering the point-to-point last passage time instead.
We prove Theorem~\ref{thm:AA} in a similar way by considering the point-to-point last passage time of a DLPP  model with a certain rotational symmetry.


A direct proof of Theorem~\ref{thm:A} was recently established in \cite{Corwin-Quastel-Remenik13}. 
The authors of \cite{Corwin-Quastel-Remenik13} first extended the determinantal formula for the joint distribution of $\Aip(\tau)$ at finitely many times to a determinantal formula for $\Prob(\Aip(\tau)\le g(\tau), \, t\in [-T,T])$ for general function $g$ and $T>0$, and then showed how Theorem~\ref{thm:A} can be obtained from this formula when $g(\tau)=\tau^2+x$. 
It is an interesting open question to find a similar direct proof for Theorem~\ref{thm:Aab} and Theorem~\ref{thm:AA} as well as Theorem~\ref{thm:AB} and Theorem~\ref{thm:ABb} below. 
The reader is referred to \cite{Quastel-Remenik13} for a survey of this direct proof and also other identities for the cousins of the Airy process. 

There is one more known identity for the Airy process. 


\begin{thm}[\cite{Quastel-Remenik12}]\label{thm:AAh}
For every $w\in \R$ and $x\in \R$, 
\begin{equation}\label{eq:AAh}
	\Prob\left[  \sup_{\tau\le w}  \Aipp(\tau) \le x - \min\{0, w\}^2 \right] = G_w^{2\to 1}(x)
\end{equation}
where $G_w^{2\to 1}(x)$ is the marginal distribution function of the process $\mathcal{A}_{2\to 1}$ introduced in \cite{Borodin-Ferrari-Sasamoto08} (see (1.7) of \cite{Quastel-Remenik12}).
\end{thm}

The distribution $G_w^{2\to 1}(x)$ interpolates $F_2$ and $F_1$: It converges to $F_2(x)$ as $w\to -\infty$ and to $F_1(4^{1/3}x)$ as $w\to +\infty$.
The paper \cite{Quastel-Remenik12} gave a direct proof of~\eqref{eq:AAh} using the method of \cite{Corwin-Quastel-Remenik13}.
In terms of DLPP models, this identity can be obtained by considering a point-to-half line last  passage time by using the result of \cite{Borodin-Ferrari-Sasamoto08}.

\bigskip

In all of the above identities, the distribution of the argmax, $\tau_{\max}$, at which the supremum is attained is also of great interest since it describes the transversal fluctuations of the associated DLPP model. 
The distribution of $\tau_{\max}$ for the identity~\eqref{eq:thmA} was computed in two recent papers \cite{Moreno-Quastel-Remenik13} and \cite{Schehr12} independently. 
The paper \cite{Moreno-Quastel-Remenik13} is mathematical and rigorous while  \cite{Schehr12} is physical. 
The density functions obtained in these papers, which look very different, are subsequently  
shown to be the same in \cite{Baik-Liechty-Schehr12}.
It is an interesting open question to find the distribution of argmax for other identities. 


\subsection{Airy process plus Brownian motions}


The indirect method for Theorem~\ref{thm:A},~\ref{thm:Aab} and~\ref{thm:AA} can also be used to prove other identities involving the Airy process and the Brownian motion if we consider DLPP models with special rows and columns. We mention two known results. 

\begin{defn}\label{defn:TW2}
For real parameters $w_+$ and $w_-$, let $F_{st}(x; w_+, w_-)$ denote the distribution function defined as 
$H(x; w_+, w_-)$ in Theorem 3.3 and Definition 3 of  \cite{Baik-Rains00}.
\end{defn}

An explicit formula is included in Section~\ref{sec:formula} below. 
The  function $F_{st}(x; w_+, w_-)$ is symmetric in the parameters $w_+$ and $w_-$. 
This function is the limiting distribution for the fluctuations of the height of totally exclusion processes starting with   Bernoulli initial conditions (see, e.g., \cite{Prahofer-Spohn02a}). 
The parameters $w_+$ and $w_-$ are associated to the initial density of the particles on the positive and negative parts,  respectively.
Especially $F_{st}(x;w , -w)$ appears when the initial condition is random and stationary. 
The following identity is obtained recently. 


\begin{thm}[\cite{Corwin-Liu-Wang13}]\label{thm:AB}
Let $\mathcal{B}(\tau)$, $\tau\in\R$, be a two-sided standard Brownian motion with $\mathcal{B}(0)=0$. Then 
\begin{equation}\label{eq:AB1}
  	\Prob\left[\sup_{\tau\in \R} \left(\Aipp(\tau)+\sqrt{2}\mathcal{B}(\tau) +4(w_+1_{\tau<0}-w_-1_{\tau>0})\tau\right)\le x\right]=F_{st}(x; w_+, w_-)
\end{equation}
for every $x\in \R$ and $w_+, w_-\in \R$. 
\end{thm}

Setting $w_+=w$, $w_-=-w$, replacing $\hat{A}(\tau)$ by $\Aip(\tau-2w)-\tau^2$ using the stationarity of the Airy process, and replacing $x$ by $x+4w^2$, we obtain 
\begin{equation}\label{eq:AB11}
  	\Prob\left[\sup_{\tau\in \R} \left(\Aip(\tau-2w)-(\tau-2w)^2+\sqrt{2}\mathcal{B}(\tau)  \right)\le x\right]
	=F_{st}(x+4w^2; w, -w).
\end{equation}
This is the one-dimensional distribution of the conjecture in the displayed equation between (1.31) and (1.32) of \cite{Quastel-Remenik13}. 
 The process $\mathcal{A}_{st}(w)$ in \cite{Quastel-Remenik13} has one-dimensional distribution function
$F_{st}(x+w^2; \frac{w}{2}, -\frac{w}{2})$.

The next result requires the following definition. 

\begin{defn}\label{defn:TW3}
For $k=1,2,\cdots$, let $F^{spiked}_k(x ; w_1, \cdots, w_k)$ denote the distribution function defined as $F_k(x; w_1, \cdots, w_k)$ in (54) of \cite{Baik-Ben_Arous-Peche05} or also equivalently in Corollary 1.3 of \cite{Baik06}.
\end{defn}

These functions are invariant under the permutations of the parameters $w_1, \cdots, w_k$.
They are the limiting distributions of the fluctuations of the largest eigenvalue of the so-called spiked random matrix models.  
It is known that $F_{st}(x; w_+, w_-)\to F^{spiked}_1(x; w_+)$ as $w_-\to +\infty$ (see (3.23) of \cite{Baik-Rains00} and (1.21) of \cite{Baik06}) and 
$F^{spiked}_1(x;0)= F_1(x)^2$ (see (24) of \cite{Baik-Ben_Arous-Peche05}).

\begin{thm}[\cite{Corwin-Liu-Wang13}]\label{thm:ABb}
Let $\mathcal{B}_1(\tau), \mathcal{B}_2(\tau), \cdots$, $\tau\ge 0$, be independent standard Brownian motions.
Let $w_1, w_2, \cdots$ be real parameters and 
let $\hat{\mathcal{B}}_i(\tau; w_i)=\mathcal{B}_i(\tau)-2\sqrt{2}w_i\tau$, the Brownian motion with drift. Then for every $x\in \R$ and $k=1,2,\cdots$, 
\begin{equation}\label{eq:spikedresult}
\begin{split}
    &\Prob\left[\sup_{0=\tau_0\le \tau_1\le \cdots\le \tau_k} \left(\Aipp(\tau_k)
    +\sqrt{2} \sum_{i=1}^k 
    \left(\hat{\mathcal{B}}_i(\tau_i; w_i)- \hat{\mathcal{B}}_i(\tau_{i-1}; w_i) \right) \right) \le x\right]
    = F^{spiked}_k (x; w_1, \cdots, w_k).
\end{split}
\end{equation}
In particular, if $\mathcal{B}(\tau)$, $\tau\ge 0$, is a standard Brownian motion, 
\begin{equation}
    \Prob\left[\sup_{\tau\ge 0} \left(\Aipp(\tau)+\sqrt{2}\mathcal{B}(\tau) \right)\le x\right]=F_1(x)^2.
\end{equation}
\end{thm}

\section{Proof of Theorem~\ref{thm:AA}}

As mentioned above, this proof is similar to that of \cite{Johansson03} for Theorem~\ref{thm:A}. 

Fix a parameter $q\in (0,1)$.
We use the notation $X\sim \Geom(q)$ to indicated that $X$ is a (shifted) geometric random variable which has the probability mass function $(1-q)q^k$, $k=0,1,2,\cdots$.

In \cite{Baik-Rains01}, five types of symmetries of DLPP models with geometrically distributed weights were considered and the limit law for the fluctuations of the last passage time was obtained for each case. 
One of the symmetry types is the rotational symmetry described as follows. 
Let $R_N:= \{(i,j): i,j= 1, \cdots, 2N\}$. 
To each $(i,j)\in R_N$ we associate a weight $w(i,j)$ with the condition that 
\begin{equation}\label{eq:symc}
	w(i,j)= w(2N+1-i, 2N+1-j), \qquad (i,j)\in R_N. 
\end{equation}
Apart from the above symmetry conditions, we assume that the random variables are independent and are distributed as $\Geom(q)$, that is, $w(i,j)$ for $(i,j)$ in $H_N:=\{(i,j) : j< 2N+1-i\}\cup\{(i, 2N+1-i): i=1, \cdots, N\}$ are independent and distributed as $\Geom(q)$, 
and $w(i,j)$ for $(i,j)\in R_N\setminus H_N$ are defined by the symmetry condition~\eqref{eq:symc}.

Let $G^\symmUU(N)$ denote the last passage time from $(1,1)$ to $(2N, 2N)$ in this model: 
\begin{equation}
	G^\symmUU(N) := \max_{\pi} W(\pi), \qquad W(\pi):= \sum_{(i,j)\in \pi} w(i,j)
\end{equation}
where $\pi$ is an up/right path from $(1,1)$ to $(2N, 2N)$ consisting of sites whose coordinates increase weakly
and $W(\pi)$ is the weight of path $\pi$. 
Here the superscript $\symmUU$ represents the symmetry of the model under the $180$ degree rotation about the center of the square. 
This model is equivalent to the model (ii)  in Section 4.4 of \cite{Baik-Rains01} (see page 14) after simple translations.
It was shown in Theorem 4.2 (ii) of \cite{Baik-Rains01} that 
\begin{equation}\label{eq:limitlaw}
	\lim_{N\to \infty} \Prob\left[  \frac{G^\symmUU(N)- (2N)\mu}{2^{2/3}(2N)^{1/3}\sigma} \le x \right] = F_2(x)^2
\end{equation}
where
\begin{equation}\label{eq:musigma}
	\mu:= \frac{2\sqrt{q}}{1-\sqrt{q}}, \qquad \sigma:=\frac{q^{1/6}(1+\sqrt{q})^{1/3}}{1-\sqrt{q}}.
\end{equation}
The Poisson version of this percolation model is related to the combinatorial problem of finding the longest increasing subsequence of random signed permutations and also the so-called $2$-colored permutation 
(see, e.g. Remark 2 in Section 1 of \cite{Baik-Rains01}), and the limit theorem analogous to the above was obtained in  \cite{Tracy-Widom99, Borodin99}.


We now consider $G^\symmUU(N)$ in a different way. 
Let $\mathcal{L}_N:= \{(i, 2N+1-i): i=1, \cdots, 2N\}$ be the sites on the straight line joining $(1, 2N)$ and $(2N, 1)$. 
Every up/right path from $(1,1)$ to $(2N, 2N)$ intersects $\mathcal{L}_N$ at a unique point. 
If we consider the paths which  intersect $\mathcal{L}_N$ at a specific point $(i, 2N+1-i)$, 
then the maximal weight among these paths equals the 
sum of two point-to-point last passage times, one from $(1,1)$ to $(i, 2N+1-i)$ and the other from $(i, 2N+1-i)$ to $(2N, 2N)$, minus the weight at $(i, 2N+1-i)$, which is counted twice. 
But due to the symmetry conditions~\eqref{eq:symc}, the last passage time from $(i, 2N+1-i)$ to $(2N, 2N)$ 
equals the last passage time from $(1,1)$ to $(2N+1-i, i)$.
Therefore, considering all possible intersection points, we find that 
$G^\symmUU(N)$ has the same distribution as 
\begin{equation}\label{eq:maxform}
	\max_{-N+1\le u\le N} \left( G(N+u, N+1-u)+G(N+1-u, N+u) \right) +D_N 
\end{equation}
where $G(m,n)$ is the usual point-to-point last passage time from $(1,1)$ to $(m,n)$ with i.i.d. (shifted) geometric weights with no symmetry conditions 
and $|D_N|$ is bounded by the maximum of $N$ independent $\Geom(q)$ random variables, which is due to the double-counting of the weights on $\mathcal{L}_N$.
Since $D_N=O((\log N)^{1+\epsilon})$ with probability $1$ as $N$ tends to infinity for any $\epsilon>0$, we obtain,  inserting~\eqref{eq:maxform} into~\eqref{eq:limitlaw}, that 
\begin{equation}\label{eq:long00}
	\lim_{N\to \infty} \Prob\left[  \frac{\max_{-N+1\le u\le N} \left( G(N+u, N+1-u)+G(N+1-u, N+u) \right)- (2N)\mu}{2^{2/3}(2N)^{1/3}\sigma} \le x \right] = F_2(x)^2.
\end{equation}
Since $G(m,n+1)$ and $G(m+1, n)$ are between $G(m,n)$ and $G(m+1, n+1)$, 
we can change $N+1$ to $N$ and find 
\begin{equation}\label{eq:long}
\begin{split}
	&\lim_{N\to \infty} \Prob\left[  \frac{\max_{|u|<N} \left( G(N+u, N-u)+G(N-u, N+u) \right)- (2N)\mu}{2^{2/3}(2N)^{1/3}\sigma} \le x \right] = F_2(x)^2.
\end{split}
\end{equation}
We show that  the left-hand side equals the left-hand side of~\eqref{eq:AA}. 

The limit of $G(m,n)$ is well known. 
It was shown in \cite{Johansson00} that 
\begin{equation}\label{eq:GLLN}
	\frac{G(m,n)}{qm+qn+2\sqrt{qmn}}  \to \frac1{1-q}
\end{equation}
in expectation and also in probability as $m,n\to \infty$ such that $m/n$ is in a compact subset of $(0, \infty)$. 
This implies, after a simple calculation, that the maximum in~\eqref{eq:maxform} is attained at $u$ near $0$. 
This statement can be made precise: see~\eqref{eq:tailestim} below. 
Now if we scale $u= N^{2/3}\tau$ and consider $G(N+ N^{2/3}\tau, N- N^{2/3}\tau)$ as a process in ``time'' $\tau$, then its fluctuations converge to the time-scaled Airy process minus a parabola in the finite distribution sense as well as in the functional limit sense \cite{Johansson03}. 
(Here $N^{2/3}$ is the scale of the so-called transversal fluctuations of KPZ universality class \cite{Johansson-transv, Baik-Deift-McLaughlin-Miller-Zhou}.)
More precisely, if we set\footnote{In (1.8) of \cite{Johansson03}, $\sigma$ is given as $\frac{q^{1/6}(1+\sqrt{q})^{1/3}}{1-q}$. This is a typographical error. The correct formula of $\sigma$ is $\frac{q^{1/6}(1+\sqrt{q})^{1/3}}{1-\sqrt{q}}$ as 
in \cite{Johansson00}.}
\begin{equation}\label{eq:Hnotation}
  	H_N(\tau):=\frac{G( N+d^{-1} N^{2/3}\tau, N- d^{-1}  N^{2/3}\tau)-\mu  N}{\sigma  N^{1/3}},
	\qquad d:= \frac{q^{1/6}}{(1+\sqrt{q})^{2/3}}, 
\end{equation}
then $H_N(\tau)$ converges to $\Aip(\tau)-\tau^2$ in the sense of weak convergence of the probability measure on $C[-T, T]$ for every fixed $T>0$. 
Substituting~\eqref{eq:Hnotation} into~\eqref{eq:long}, we obtain
\begin{equation}\label{eq:zl-7-9-1}
	\lim_{N\to \infty} \Prob\left[  \max_{|\tau|< dN^{1/3}} \left( \frac{H_N(\tau)+H_N(-\tau)}{2} \right) \le x \right] = F_2(x)^2.
\end{equation}
It remains to show that the left-hand side equals the left-hand side of~\eqref{eq:AA}. 

The functional limit theorem mentioned above implies that  
\begin{equation}\label{eq:functional}
	\lim_{N\to\infty}\Prob\left[ \max_{|\tau|\le T} \left( \frac{H_N(\tau)+H_N(-\tau)}{2} \right) \le x\right]
=\Prob\left[ \max_{|\tau|\le T} \left( \frac{\Aipp(\tau)+\Aipp(-\tau)}{2} \right) \le x\right]
\end{equation} 
for any fixed $T>0$.
The tail estimate of $H_N(\tau)$ for large $|\tau|$ was also obtained in \cite{Johansson03}. 
The analysis in that paper implies that (see (127) in \cite{Baik-Liu02} for details) for every fixed $x\in\R$ and $\epsilon>0$, 
there are $T_0$ and $N_0$ such that 
\begin{equation}\label{eq:tailestim}
	\Prob\left[\max_{|\tau|>T}H_N(\tau)>x\right]<\epsilon
\end{equation}
for all $T>T_0$ and $N>N_0$. 
Since $\Prob\left[\max_{|\tau|\le T} \left( \frac{\Aipp(\tau)+\Aipp(-\tau)}{2} \right) \le x\right]$ is monotone in $T$,~\eqref{eq:functional} and~\eqref{eq:tailestim} imply that 
the left-hand side of~\eqref{eq:zl-7-9-1} equals the left-hand side of~\eqref{eq:AA}. 
Theorem~\ref{thm:AA} is proved.

\begin{rmk}
In addition to the symmetry $\symmUU$, four other symmetry types were considered in \cite{Baik-Rains01}.
They are indicated by symbols ${\symmU}$, $\symmS$, $\symmO$, and $\symmu$.
The first one $\symmU$ has no symmetry. 
If we consider the last passage time in this case as above, we arrive at Theorem~\ref{thm:Aab}: the maximal path in $T_N:=\{(i,j): i+j\le N+1\}$ and $R_N\setminus T_N$ give rise to two independent Airy processes since the last passage times in two parts are independent to each other. 
The second symmetry type $\symmS$ has the reflection symmetry about the line $\mathcal{L}_N$. Since the maximal path in $T_N$ and $R_N\setminus T_N$ are identical except for the asymptotically negligible contribution from the site on $\mathcal{L}_N$, their sum is basically twice of the maximal path in $T_N$, thus giving rise to single Airy process. This leads to Theorem~\ref{thm:A}. 
\end{rmk}

\begin{rmk}
The symmetry types $\symmO$, and $\symmu$ contain an extra reflection symmetry about the other diagonal line $\{(i,j): j=i\}$. 
Hence we can consider these cases as directed last passage models in the triangle $\mathcal{T}:= \{(i,j): j\le i\}$. 
The limiting distribution for $\symmO$ and $\symmu$ is $F_4$ and $F_2$, respectively. 
Now the analog of~\eqref{eq:Hnotation} in triangle $\mathcal{T}$ does not converge to Airy process but to a different process which interpolate $F_4$ and $F_2$ \cite{Imamura-Sasamato04a}.
Thus we can expect that the the supremum of two independent such processes minus a parabola, over $\tau\in [0, \infty)$, equals $F_4$, and the supremum of one such process minus parabola equals $F_2$. 
To make this statement rigorous, we need the functional theorem and the tail estimate for DLPP in triangle $\mathcal{T}$. This will be discussed in future work. 
We note that by making the weights on the diagonal line $\{(i,j): j=i\}$ different from the rest of the triangle, 
we can obtain identities for yet another process.
\end{rmk}

\section{Outline of proof of Theorem~\ref{thm:AB} and~\ref{thm:ABb}}

In this section, we outline the proof of Theorem~\ref{thm:AB} and~\ref{thm:ABb} obtained in \cite{Corwin-Liu-Wang13} in order to illustrate how different DLPP models give rise to different identities. 
The basic idea is same as the previous section and \cite{Johansson03}: we consider a DLPP model for which the limit theorem is proved, and then interpret the last passage time as the maximum of last passage times of paths with arbitrary end-points. 
The technical part is the functional limit 
and the tail estimate. 
It turned out that for Theorem~\ref{thm:AB} and~\ref{thm:ABb} one needs the functional limit theorem and the tail estimates along horizontal lines and vertical lines instead of the diagonal line in the previous section. 
Such results are obtained in \cite{Corwin-Liu-Wang13} by proving the so-called slow decorrelation phenomenon 
and using the results of Johansson \cite{Johansson03}.

For Theorem~\ref{thm:AB}, one uses the following DLPP model: the weights are independent and satisfy
\begin{equation}
\begin{split}
	&w(i,j)\sim \Geom(q), \qquad\qquad   i,j=1, \cdots, N, \\
	&w(i,0)\sim \Geom(\alpha_+ \sqrt{q}), \qquad i=1,\cdots, N, \\
	&w(0,j)\sim \Geom(\alpha_-\sqrt{q}), \qquad j=1,\cdots,  N, \\
	&w(0,0)=0
\end{split}
\end{equation}
where $q\in (0,1)$ is a fixed parameter and 
\begin{equation}\label{eq:aplha}
\begin{split}
	\alpha_+= 1- \frac{2w_+}{\sigma N^{1/3}},  \qquad 
	\alpha_-= 1- \frac{2w_-}{\sigma N^{1/3}}, 
\end{split}
\end{equation}
for fixed real parameters $w_+$ and $w_-$. Here $\sigma$ is same as~\eqref{eq:musigma}. 
It was shown in Section 4 of \cite{Baik-Rains00} that 
the last passage time $X(N)$ from $(0,0)$ to $(N,N)$ satisfies 
\begin{equation}\label{eq:limitst}
\begin{split}
	\lim_{N\to \infty} \Prob\bigg[ \frac{X(N)-\mu N}{\sigma N^{1/3}} \le x \bigg] 
	= F_{st}(x; w_+, w_-).
\end{split}
\end{equation}

We now consider the last passage time in a different way. The last passage path, considered as starting at $(N,N)$ and ending at $(0,0)$, either arrives at a point $(i,0)$ for some $i$ and travel horizontally left to $(0,0)$, or arrives at a point $(0,j)$ for some $j$ and travel vertically downward to $(0,0)$. 
Hence, denoting by $G(N-i, N-j)$ the last passage time from $(N,N)$ to $(i,j)$, 
we find that 
\begin{equation}\label{eq:Xmax}
\begin{split}
	X(N) = \max \left\{  \max_{i=1, \cdots, N} \left(G(N-i, N)+ S_i^+ \right), 
	\max_{j=1, \cdots, N}  \left(G(N, N-j)+ S_j^-\right) \right\} 
\end{split}
\end{equation}
where $S_i^+$ is the sum of $i$ independent $\Geom(\alpha_+\sqrt{q})$ random variables and $S_j^-$ is the sum of $j$ independent $\Geom(\alpha_-\sqrt{q})$ random variables.
From~\eqref{eq:GLLN} and the law of large numbers of independent variables, it is reasonable to expect that the maximum of the above expression occurs when $i$ and $j$ are close to $0$. 
Now set 
\begin{equation}
\begin{split}
	\tilde{H}_N(\tau):= \frac{G(N-2d^{-1}N^{2/3}\tau, N) -\mu (N -  d^{-1} N^{2/3} \tau )}{\sigma N^{1/3}}, \qquad \tau\ge 0, \\
	\tilde{H}_N(\tau):= \frac{G(N, N+2d^{-1}N^{2/3}\tau) -\mu (N +  d^{-1} N^{2/3} \tau) }{\sigma N^{1/3}}, \qquad \tau\le 0, \\
\end{split}
\end{equation}
where $d$ is defined in~\eqref{eq:Hnotation}.
Then $\tilde{H}_N(\tau)$ converges to $\Aipp(-\tau)$ in the sense of finite distribution and also in the functional limit sense \cite{Corwin-Liu-Wang13}. 
On the other hand, since the mean and the variance of $\Geom(a)$ is $\frac{a}{1-a}$ and $\frac{a}{(1-a)^2}$ respectively, we find from the central limit theorem, after inserting~\eqref{eq:aplha}, that 
\begin{equation}
\begin{split}
	\frac{S_{2d^{-1}N^{2/3}\tau} ^{\pm}- \mu d^{-1} N^{2/3} \tau}{ \sigma N^{1/3}} 
	\Rightarrow  \sqrt{2} \mathcal{B}_{\pm}(\tau) - 4 w_{\pm} \tau
\end{split}
\end{equation}
in distribution where $\mathcal{B}_+(\tau)$ and $\mathcal{B}_-(\tau)$ are independent standard Brownian motions. 
Thus, we find that, at least formally, $\frac{X(N)-\mu N}{\sigma N^{1/3}}$ converges to 
\begin{equation}\label{eq:Xmax2}
\begin{split}
	\max \left\{  \sup_{\tau\ge 0} \left( \Aipp(-\tau) +  \sqrt{2} \mathcal{B}_+(\tau)- 4 w_+ \tau \right),  
	\sup_{\tau\le 0} \left(\Aipp(-\tau) +  \sqrt{2} \mathcal{B}_-(-\tau)+ 4 w_- \tau \right) \right\} .
\end{split}
\end{equation}
This argument was made rigorous in \cite{Corwin-Liu-Wang13}.
After we change $\tau$ to $-\tau$, this, combined with~\eqref{eq:limitst}, implies Theorem~\ref{thm:AB}.

\bigskip

Theorem~\ref{thm:ABb} is obtained by considering the DLPP model where the weights are independent and satisfy
\begin{equation}
\begin{split}
	&w(i,j)\sim \Geom(q), \qquad\qquad   i=k+1, \cdots, N, \quad ,j=1, \cdots, N, \\
	&w(i,j)\sim \Geom(\alpha_i\sqrt{q}), \qquad i=1, \cdots,  k, \quad j=1,\cdots,  N
\end{split}
\end{equation}
for fixed integer $k$ where
\begin{equation}\label{eq:aplha23}
\begin{split}
	\alpha_i= 1- \frac{2w_i}{\sigma N^{1/3}},  \qquad i=1,\cdots, k. 
\end{split}
\end{equation}
The limit theorem to $F^{spiked}(x; w_1, \cdots, w_k)$, similar to~\eqref{eq:limitst}, was proved in \cite{Baik-Ben_Arous-Peche05} for exponentially distributed weights and in Theorem 2-3$'$ \cite{Imamura-Sasamoto07} (set $m=1$) for geometrically distributed weights. 
On the other hand, the analog of~\eqref{eq:Xmax} is 
\begin{equation}\label{eq:Xmax00}
\begin{split}
	 \max_{0=j_0\le j_1\le \cdots \le j_k\le N}  \left\{  \sum_{i=1}^k  (S_{j_{i}}^{(i)} - S_{j_{i-1}-1}^{(i)} )
	 + G(N-k, N-j_k) \right\} 
\end{split}
\end{equation}
where $S_{\ell}^{(i)}= w(i, 1)+ \cdots +w(i, \ell)$. 
This leads to the left-hand side of~\eqref{eq:spikedresult}. See \cite{Corwin-Liu-Wang13} for the detail.

\section{Formula of distribution functions}\label{sec:formula}

For the convenience of the reader, we include the formulas of the Tracy-Widom distribution functions 
and also the distribution functions in Definition~\ref{defn:TW2} and~\ref{defn:TW3}. 
All of them have at least two different expressions. One expression involves the Painlev\'e equation and another expression involves a Fredholm determinant of an operator whose kernel is related to the Airy function. 
Here we only present the formulas involving the Painlev\'e equation. 

Let $q(x)$ be the  solution to the Painlev\'e II equation 
$q''=2q^3+xq$ satisfying the condition that $q(x) \sim \Ai(x)$ as $x\to +\infty$ where $\Ai(x)$ is the Airy function.\footnote{The function $u(x)=-q(x)$ is used in \cite{Baik-Rains00} and \cite{Baik06}} 
The solution is unique and smooth, and is called the  Hastings-McLeod solution 
\cite{HastingsMcLeod, Fokas-Its-Kapaev-Novok}. 
The Tracy-Widom distributions are defined as \cite{Tracy-Widom94, Tracy-Widom96}
\begin{equation}
\begin{split}
	F_2(x)= F(x)^2, \qquad F_1(x)= F(x)E(x), \qquad F_4 \left(x/\sqrt{2} \right)= \frac12 F(x) \left( E(x)+\frac1{E(x)} \right)
\end{split}
\end{equation}
where
\begin{equation}
\begin{split}
	F(x)= e^{-\frac12\int_x^\infty (s-x) q(s)^2 ds}, \qquad 
	E(x)= e^{-\frac12 \int_x^\infty q(s)ds}.
\end{split}
\end{equation}

The distribution functions in Definition~\ref{defn:TW2} and~\ref{defn:TW3} are more involved. 
Nevertheless they are expressible only in terms of $q(x)$ above and two other functions $a(x;w)$ and $b(x;w)$.
Let $a(x;w)$ and $b(x;w)$ be the solution to the initial value problem of the system of first order linear differential equations, 
\begin{equation}
\begin{split}
	\frac{d}{dx} \begin{pmatrix} a(x;w) \\ b(x;w) \end{pmatrix} = \begin{pmatrix}
	0 & -q(x) \\ -q(x) & -2w \end{pmatrix} \begin{pmatrix} a(x;w) \\ b(x;w) \end{pmatrix}, 
	\qquad 
	\begin{pmatrix} a(0;w) \\ b(0;w) \end{pmatrix} = \begin{pmatrix} E(x)^2 \\ -E(x)^2 \end{pmatrix}.
\end{split}
\end{equation}
There is a unique solution which is smooth in $(x,w)\in \R\times \mathbb{C}$. 
The above differential equations are the first part of the Lax pair for the Painlev\'e II equation in the theory of integrable systems. 
Then (see (3.22) of \cite{Baik-Rains00})
\begin{equation}
\begin{split}
	F_{st}(x;w_+, w_-)	
	= F_2(x) \left( a(x; w_+)a(x; w_-) - \frac{a(x; w_+)a(x; w_-) - b(x; w_+)b(x; w_-) }{2(w_++w_-)}
	p(x) \right)
\end{split}
\end{equation}
where
\begin{equation}
\begin{split}
	p(x):= \int_\infty^x q(y)^2dy = q(x)^4+xq(x)^2-(q'(x))^2 .
\end{split}
\end{equation}
When $w_++w_-=0$, the above can also be written as  (see (3.35) of \cite{Baik-Rains00})
\begin{equation}
\begin{split}
	F_{st}(x;w, -w)	
	= \frac{d}{dx}  \left( F_2(x)\int_{-\infty}^x a(y; w)a(y; -w) dy \right).
\end{split}
\end{equation}
We note that the expectation of $F_{st}(x+4w^2; w, -w)$ is $0$ (see Proposition 3.4 in \cite{Baik-Rains00}). 
When $w_+=w_-=0$, we have (see (2.16) of \cite{Baik-Rains00}) 
\begin{equation}
\begin{split}
	F_{st}(x; 0,0)= \left\{ 1 +  (x-2q'(x)+2(q(x))^2) \left( \int_x^{\infty} q(s)^2 ds \right) \right\} (E(x))^4 F_2(x).
\end{split}
\end{equation}

The Painlev\'e formula for the distribution functions in Definition~\ref{defn:TW3} was obtained in \cite{Baik06}. 
We have (see Corollary 1.3, Lemma 1.4, and and (1.16) of \cite{Baik06})
\begin{equation}
\begin{split}
	F^{spiked}_k(x;w_1, \cdots, w_k) 
	= F_2(x) \frac{\det\left( (w_i+ \frac{d}{dx})^{j-1} f(x; w_i) \right)_{i,j=1}^k }{\prod_{1\le i<j\le k} (w_j-w_i)}, 
	\qquad f(x;w):= a(x; w/2).
\end{split}
\end{equation}
Some special cases are (see (1.40) of \cite{Baik06}) 
\begin{equation}
\begin{split}
	F^{spiked}_1(x;0) &= F_2(x) (E(x))^2 = (F_1(x))^2, \\
	F^{spiked}_{2}(x; 0,0) &=  F_2(x) (E(x))^4  \left( 1- q(x+2q^2-2q')   \right), \\
	F^{spiked}_3(x;0,0,0) &= F_2(x)(E(x))^6 \left( 1- 2q(x+2q^2-2q')  +\frac12 (q^2+q')(x+2q^2-2q')^2\right).
\end{split}
\end{equation}

\subsubsection*{Acknowledgments}
The work of Jinho Baik was supported in part by NSF grants DMS1068646.

\def\cydot{\leavevmode\raise.4ex\hbox{.}}


\end{document}